\documentclass{article}      
\usepackage{amssymb}
\usepackage{amsmath}
\usepackage{amscd}
\usepackage{amsthm}  
\usepackage{epsfig}  
\usepackage{subfigure}
\usepackage[all]{xy}
\usepackage{texdraw}
\SelectTips{cm}{}
\usepackage{wrapfig}
\newtheorem{theorem}{Theorem}
\newtheorem{corollary}[theorem]{Corollary}
\newtheorem{lemma}[theorem]{Lemma}  
\newtheorem{proposition}[theorem]{Proposition}

\theoremstyle{definition}
\newtheorem{definition}{Definition}
\newtheorem*{ques}{Question}
\newtheorem*{convent}{Conventions}

\newcommand{\br}{\mathbb{R}}

\newcommand{\p}{\partial}

\newcommand{\cl}{\mathcal{L}}

\newcommand{\cf}{\mathcal{F}}

\newcommand{\ct}{\mathcal{T}}

\newcommand{\hk}{\hookrightarrow}

\newcommand{\med}{\medskip}
\newcommand{\la}{\longrightarrow}
\newcommand{\bfl}{\begin{flushleft}}
\newcommand{\efl}{\end{flushleft}}

\newcommand{\eps}{\epsilon}
\newcommand{\xr}{\xrightarrow}

\newcommand{\lmm}{LM \times_M LM}
\newcommand{\fmk}{F_{D^k}(M, 2)}
\newcommand{\fmak}{F_{D^k}(M_1, 2)}
\newcommand{\fmbk}{F_{D^k}(M_2, 2)}
\def\:{\colon}

\newcommand{\bedka}{D(\tau_{M_1} \oplus \eps^k) }
 
\newcommand{\bedkb}{D(\tau_{M_2} \oplus \eps^k) }

\newcommand{\mkb}{M(k)_\bullet}
\newcommand{\mokb}{M_1(k)_\bullet}
\newcommand{\mtkb}{M_2(k)_\bullet}

\begin{document}  

  \title{The homotopy 
invariance of the string topology loop product and string bracket}
  \author{Ralph L. Cohen \quad John Klein \quad Dennis Sullivan \thanks{All three authors were partially supported by  grants from the NSF}}
\date{\today}
\maketitle

 
 \begin{abstract}
Let $M^n$ be a  closed, oriented,  $n$-manifold, and $LM$ its free loop space.
  In \cite{chassullivan}  a commutative algebra structure in homology, $H_*(LM)$, and a Lie algebra structure in equivariant homology 
$H^{S^1}_*(LM)$, were defined.  In this paper we prove that these structures are homotopy invariants   in the following sense.   
  Let  $f : M_1 \to M_2$ be a homotopy equivalence of closed, oriented $n$-manifolds. Then the induced equivalence, $Lf : LM_1 \to LM_2$ induces a ring isomorphism in homology, and an isomorphism of Lie algebras in equivariant homology.  The analogous statement also holds true for any generalized homology theory $h_*$ that supports an orientation of the $M_i$'s.
\end{abstract}
 \section*{Introduction}

The term ``string topology" refers to multiplicative structures on the (generalized) homology of spaces of paths and loops in a manifold.  Let $M^n$ be a closed, oriented, smooth   $n$-manifold.   The basic ``loop homology algebra" is defined by a product
$$
\mu : H_*(LM) \otimes H_*(LM) \la H_*(LM)
$$
of degree $-n$, and 
and  the ``string Lie algebra"  structure is  defined by  a bracket
$$
[\,\, , \,\, ] : H^{S^1}_*(LM) \otimes  H^{S^1}_*(LM) \la  H^{S^1}_*(LM)
$$
of degree $2-n$.  These  were   defined in \cite{chassullivan}.  Here  $H^{S^1}_*(LM)$
 refers to the equivariant homology, $ H^{S^1}_*(LM) = H_*(ES^1 \times_{S^1} LM)$.  
More basic structures on the chain level were also studied in \cite{chassullivan}.  
Furthermore, these structures were shown to exist  for any multiplicative homology 
theory $h_*$ that supports an orientation of $M$. (see \cite{cohenjones}.
 \footnote{The string bracket for generalized homology theories was not explicitly discussed in \cite{cohenjones}, although   in Theorem 2 of that paper  there is a homotopy theoretic action of Voronov's cactus operad given, which, by a result of Getzler \cite{getzler}  yields a Batalin-Vilkovisky structure on the generalized homology, $h_*(LM)$, when $M$ is $h_*$-oriented. According to \cite{chassullivan}, this is all that is needed to construct the string bracket.  We will review this construction below.})
 Alternative descriptions
of the basic structure were given in \cite{cohenjones} and \cite{chataur}, 
but in the end they all relied on various perspectives of intersection theory of   
chains and homology classes.  

The existence of various descriptions of these operations
leads to the following:

\begin{ques} {\sl To what extent are the
the string topology operations 
sensitive to the smooth structure of the manifold, 
or even the homeomorphism
structure?}
\end{ques}

The main goal of this paper is to settle this question in the case of two of the basic operations:   the string topology loop product and string bracket.  We will in fact prove more:   we will show that 
the loop homology algebra and string Lie algebra structures are oriented
{\it homotopy invariants.}  

We remark that it is still not known whether the full range of string topology operations \cite{chassullivan}, \cite{chassullivan2}, \cite{sullivan}, \cite{cohengodin} are homotopy invariants.  Indeed the third author has conjectured that they are not (see the postscript to \cite{ranickisullivan}).  More about this point will be made in the remark after theorem 2 below.

\med To state the main result, let $h_*$ be a multiplicative 
homology theory that supports an orientation of $M$. Being a multiplicative theory 
  means that the corresponding cohomology theory, $h^*$, 
admits a cup product, or  more precisely, 
the representing spectrum of the theory  is required to be a ring spectrum.  
An {\it $h_*$-orientation} of a closed $n$-manifold $M$ 
can be viewed as a choice of fundamental class  $[M] \in h_n(M)$ that
induces a Poincar\'e duality isomorphism.

\med
\begin{theorem}\label{main}  Let $M_1$ and $M_2$ be closed, $h_*$-oriented $n$-manifolds.  Let  $f : M_1 \to M_2$ be an $h_*$-orientation preserving  homotopy equivalence.       Then the induced homotopy equivalence of loop spaces, 
$Lf : LM_1 \to LM_2$ induces a ring isomorphism of loop homology algebras,
$$
\begin{CD}
(Lf)_*: h_*(LM_1) @>\cong >> h_*(LM_2).
\end{CD}
$$
Indeed it is an isomorphism of Batalin-Vilkovisky ($BV$) algebras.  Moreover the induced map on equivariant homology,
$$
\begin{CD}
(Lf)_*: h_*^{S^1}(LM) @>\cong >> h_*^{S^1}(LM).
\end{CD}
$$
is an isomorphism of graded Lie algebras,
\end{theorem}

\med
Evidence for the above theorem came from the 
results of \cite{cohenjones} and \cite{cohen} which said
that for simply connected manifolds  $M$,  
there is an  isomorphism of graded 
algebras, 
$$
H_*(LM) \cong H^*(C^*(M), C^*(M)) \, ,
$$
where the right hand side is the Hochschild cohomology of $C^*(M), $
the differential graded algebra of singular cochains on $M$, with
multiplication given by cup product.   
The Hochschild cohomology algebra is clearly a homotopy invariant.

However, the above isomorphism is defined in terms of the 
Pontrjagin-Thom construction arising from the diagonal 
embedding $M^S \subset M^T$ associated with each surjection of 
finite sets $T \to S$.  Consequently, since the Pontrjagin-Thom construction uses the smooth structure,  this isomorphism {\it a priori} 
seems to be sensitive to the smooth structure.
Without an additional argument, one can only conclude from this
isomorphism that the loop homology algebra of two homotopy equivalent 
simply connected closed manifolds 
are {\it abstractly} isomorphic. In summary, to prove 
homotopy invariance in the sense of Theorem \ref{main}, 
one needs a different argument. 
\med

The argument we present here  does not need the 
simple connectivity hypothesis. This should prove of particular  interest in
the case of surfaces and $3$-manifolds.   Our argument
uses the description of the loop product $\mu$ in terms of  a 
Pontrjagin-Thom collapse map of an embedding 
$$
LM \times_M LM \hk LM \times LM
$$ 
given in \cite{cohenjones}.  
Here $\lmm$ is the subspace of $LM \times LM$ 
consisting of those pairs of loops $(\alpha, \beta)$, with $\alpha (0) = \beta (0)$. In this description  we are thinking of the loop space as the space of piecewise smooth maps $[0,1] \to M$ whose values at $0$ and $1$ agree.  This is a smooth, infinite dimensional manifold.  The differential topology of such manifolds is discussed
in \cite{cohenstacey} and \cite{chataur}.  

This description quickly reduces the proof of the theorem to the question of 
whether the homotopy type of the complement
of this embedding, $(LM \times LM) - (\lmm)$ is a stable homotopy invariant
when considered as ``a space over'' $LM \times LM$. 
By using certain pullback properties, the latter question is then further reduced 
to the question of whether the complement of the diagonal embedding,
$\Delta : M \to M \times M$, or somewhat weaker, 
the complement of the embedding
\begin{align}
\Delta_k : M &\to M \times M \times D^k \notag \\
x &\to (x, x, 0) \notag
\end{align}
is a homotopy invariant when considered as a space over $M \times M$.   
For this we develop the notion of relative smooth and Poincare embeddings.  
This is related to  the classical theory of Poincare embeddings 
initiated by Levitt \cite{Levitt_poincare} and Wall \cite{Wall},
and further developed by the second author in \cite{Klein_haef} and 
\cite{Klein_haef2}. 
However, for our purposes, the results we need can be proved directly by elementary arguments.  The results in Section 2 on relative embeddings
are rather fundamental, but don't appear in the literature.  These results may be of independent interest, and furthermore, by proving them here, we make the paper self contained.



\med

Early on in our investigation  of this topic, 
our methods led us to advertise the following question, 
which is of  interest independent of its applications to string topology.

\med 
Let $F(M,q)$ be the configuration  space of $q$-distinct, 
ordered points in a closed manifold $M$.  

\med
\begin{ques} {\sl Assume that  $M_1$ and $M_2$ be homotopy equivalent, simply connected closed $n$-manifolds.  Are $F(M_1, q)$ and $F(M_2, q)$ 
homotopy equivalent?}
\end{ques}
\med

One knows that these configuration spaces have isomorphic cohomologies 
(\cite{milgbodig}), stable homotopy types (\cite{Aouina-Klein}, 
\cite{fcohen?}) 
and  have homotopy equivalent loop spaces (\cite{cohengitler}, 
\cite{Levitt_arcs}). 
But the homotopy invariance of the configuration spaces themselves 
is not yet fully understood. For example, when $q =2$ and the manifolds
are $2$-connected, then one does have homotopy invariance
(\cite{Levitt_arcs}, \cite{Aouina-Klein}). On the other hand,
the simple connectivity assumption in the above question
is a necessity: a recent result of Longoni and Salvatore \cite{LS} shows that
for the homotopy equivalent lens spaces $L(7,1)$ and $L(7,2)$,
the configuration spaces $F(L(7,1),2)$ and
$F(L(7,2),2)$ have distinct homotopy types. 
\med

This paper is organized as follows.  In Section 1 we will reduce the proof of the main theorem to a question
about the homotopy invariance of the  complement of the diagonal embedding, 
$\Delta_k : M \to M \times M \times D^k$.  In Section 2 we  
develop the theory of relative smooth and Poincare embeddings, 
and then apply it to prove the homotopy invariance of these 
configuration spaces, and complete the proof of 
  Theorem \ref{main}.

\med

\bfl  \bf  Remark. \rm After the results of this paper were announced,  two independent proofs of the homotopy invariance of the loop homology product were found by  by Crabb \cite{crabb}, and by  Gruher-Salvatore \cite{gruhersalvatore}.
\efl

\begin{convent} 
A finitely dominated pair of spaces $(X,\partial X)$  is 
a {\it Poincare pair} of dimension $d$ if there
exists a pair $({\mathcal L},[X])$ consisting 
of a rank one abelian local coefficient
system $\mathcal L$ on $X$ and a ``fundamental class'' 
$[X] \in H_d(X,\partial X;{\mathcal L})$
such that the cap product homomorphisms
$$
\cap [X]\: H^*(X;{\mathcal M}) \to  
H_{d-*}(X,\partial X;{\cal L}\otimes {\mathcal M}) 
$$
and
$$
\cap [\partial X]\: H^*(\partial X;{\mathcal M}) \to  H_{d-1-*}(\partial X;{\cal L}\otimes {\mathcal M}) 
$$
are isomorphisms for all local coefficient bundles ${\mathcal M}$ on $X$ (respectively
on $\partial X$). Here $[\partial X]\in H_{d-1}(\partial X;{\cal L})$ denotes
the image of $[X]$ under the evident boundary homomorphism.
If such a pair  $({\mathcal L},[X])$   exists, then it is unique up to unique 
isomorphism.
\end{convent}

\med

\section{A question about configuration spaces}

\med
In this section we state one of our main results about the homotopy invariance of certain configuration spaces,  and then use it to prove Theorem \ref{main}. The theorem about configuration spaces will be proved in section 2.   

\med
Using an identification of the tangent bundle $\tau_{M } $ with the normal bundle of the diagonal, $\Delta \: M  \to M  \times M $,   we have an embedding of the disk bundle,
$$
D(\tau_{M }) \subset M \times M \, ,
$$
which is identified with a compact tubular neighborhood of
the diagonal.    (To define the unit disk bundle, we use a fixed Euclidean structure on $\tau_M$.)  
The closure of its complement will be denoted $F(M, 2)$.   Notice  that
the inclusion $F(M, 2) \subset M \times M - \Delta$ is a weak
equivalence.   We therefore have a decomposition, 
$$
M \times M =  D(\tau_{M }) \cup_{S(\tau_{M })} F(M, 2),
$$
where $S(\tau_M) = \p D(\tau_M)$ is the unit sphere bundle. 
 We now vary the configuration space in the following way.

Let $D^k$ be a  closed  unit disk, and consider the generalized diagonal embedding,
  \begin{align}
\Delta_k : M &\to M \times M \times D^k \notag \\
x &\to (x, x, 0). \notag
\end{align}

We may now identify the stabilized tangent bundle, $ \tau_M \oplus \eps^k $ with the normal bundle of this embedding, where $\eps^k$ is the trivial $k$-dimensional bundle.
 This yields an embedding, $D(\tau_M \oplus \eps^k) \subset M \times M \times D^k$, which is identified with a closed tubular neighborhood of $\Delta_k$.  The closure of its complement is denoted by $\fmk$.  The reader will notice that this is a model for the $k$-fold fiberwise suspension of the map $F(M,2) \to M \times M$.   We  now have a similar decomposition,
 $$
 M \times M \times D^k = D(\tau_M \oplus \eps^k) \cup_{S(\tau_M \oplus \eps^k)} \fmk.
 $$
 Notice furthermore, that the boundary, $\p ( M \times M \times D^k) = M \times M \times S^{k-1}$  lies in the subspace, $\fmk$.   In other words we have a commutative
 diagram,
 \begin{equation}\label{m2k}
 \begin{CD}
 S(\tau_M \oplus \eps^k)   @>>> \fmk  @<<< M \times M \times S^{k-1} \\
 @VVV   @VVV \\
 D(\tau_M \oplus \eps^k)  @>>>  M \times M \times D^k
 \end{CD}
 \end{equation}
 where the commutative square is a pushout square.  We refer to this diagram as $\mkb$.  We think of this more functorially as follows.
 
 Consider the  partially ordered set  $\cf$, with 
five objects, referred to as 
$\emptyset$, $0$,  $1$, $01$, and $b$, and the morphisms  
are generated by the following commutative diagram
 \begin{equation}\label{five}
 \begin{CD}
\emptyset  @>>> 1  @<<< b \\
@VVV    @VVV \\
0 @>>> 01\, .
\end{CD}
\end{equation}
Notice that $01$ is a terminal object of this category. 

\med
\begin{definition}\label{fivespace}
We define an $\cf$-space to be a functor  $X : \cf \to \text{\rm Top}$,
where \text{Top} is the category of 
topological spaces. The value of the functor at $S \subset \{0,1\}$ is
denoted $X_S$. It will sometimes be convenient to specify $X$
by maps of pairs
$$
(X_b,\emptyset) \to (X_1,X_\emptyset) \to (X_{01},X_0) \, ,
$$
where we are abusing notation slightly since the maps
$X_\emptyset \to X_1$ and $X_0 \to X_{01}$ need not be
inclusions.

A map (morphism) $\phi : X \to Y$ of $\cf$-spaces is a 
natural transformation of functors.  
We  say that $\phi$ is a weak equivalence, if it 
is an object-wise weak homotopy equivalence, i.e, it gives a
a weak homotopy equivalence $\phi_i :X_i \xr{\simeq}Y_i$ 
for each object $i \in \cf$.  
In general, we say that two  $\cf$-spaces 
are weakly equivalent if there is a 
finite zig-zag of morphisms connecting them, 
$$
X = X^1 \xr{} X^2 \xleftarrow{ } X^3 \xr{ }   \cdots \leftarrow  X^i \rightarrow \cdots  X^n = Y \, ,
$$
where each morphism is a  weak equivalence.   \end{definition}

 Notice that   diagram (\ref{m2k}) defines a $\cf$-space for each  closed manifold $M$, and integer $k$.  We call this $\cf$- space $\mkb$.    
In particular, $M(k)_{01} = M \times M \times D^k$. 
  
  \med
 The following is our main result about configuration spaces.  It will be proved in section 2.

\med
\begin{theorem}\label{config}  Assume $M_1$ and $M_2$ are closed
  manifolds and that
$f : M_1 \to M_2$ is a  homotopy 
equivalence.   Then for $k$ sufficiently large,  the $\cf$-spaces $M_1(k)$ and $M_2(k)$ are weakly equivalent in the following specific way.

There is a $\cf$-space $\ct_\bullet$ that takes values in spaces of the homotopy type of $CW$-complexes, and morphisms of $\cf$-spaces,
$$
\mokb \xr{\phi_1} \ct_\bullet \xleftarrow{\phi_2} \mtkb
$$
satisfying the following properties:
\begin{enumerate}
\item The morphisms $\phi_1$ and $\phi_2$ are weak equivalences.
\item The terminal space $\ct_{01}$ is defined as
$$
\ct_{01} = T_{f\times f} \times D^k
$$
where $T_{f \times f}$ is the mapping cylinder 
$(M_2 \times M_2) \cup_{f\times f}(M_1 \times M_1) \times I$.  
Furthermore on the terminal spaces, the morphisms, $\phi_1 : M_1 \times M_1 \times D^k \to T_{f\times f} \times D^k$ and $\phi_1 : M_2 \times M_2 \times D^k \to T_{f\times f} \times D^k$ are given by $\iota_1 \times 1$ and $\iota_2  \times 1$, where for $j = 1,2$, $\iota_j : M_j \times M_j \to T_{f \times f}$ are the obvious inclusions as the two ends of the mapping cylinder.
\item The induced weak equivalence,
$$
\begin{CD}
 D(\tau_{M_1} \oplus \eps^k)  
 = M_1(k)_0 @>\phi_2 > \simeq > \ct_0
  @<\phi_2<\simeq < M_2(k)_0 =  D(\tau_{M_2} \oplus \eps^k) 
 \end{CD}
 $$
 is homotopic to the composition
$$\begin{CD}\bedka  @>\text{\rm project}>\simeq > M_1 @>f>\simeq > M_2  
@>\text{\rm zero}>\simeq > \bedkb. \end{CD}$$
\end{enumerate}
\end{theorem}

 \med
 Notice that this theorem is a strengthening of the following homotopy invariance statement (see \cite{Aouina-Klein}). 
 
 \begin{corollary} Let $f : M_1 \to M_2$ be a homotopy equivalence of closed manifolds. Then for sufficiently large $k$, the configuration spaces, $\fmak$ and $\fmbk$ are homotopy equivalent.
 \end{corollary}
 
 \med
 As mentioned, we will delay the proof of Theorem \ref{config} until the next section.  Throughout the rest of this section we will assume its validity, and will use it to prove
 Theorem \ref{main}, as stated in the introduction.
 
 \begin{proof}
 Consider the equivalences of $\cf$-spaces given in Theorem \ref{config}.  Notice that we have the following commutative diagram of maps of pairs.
 \begin{equation}\label{downstairs}
 \begin{CD}
 (M_1(k)_{01}, \,  M_0(k)_b)  @>>>   
(M_1(k)_{01}, \,  M_1(k)_1)  @<<<   (M_1(k)_0 , \,  M_1(k)_\emptyset) \\
 @V\phi_1 VV  @VV\phi_1 V  @VV\phi_1 V \\
 (\ct_{01}, \,  \ct_b)  @>>>   (\ct_{01}, \,  \ct_1)  
@<<<   (\ct_0 , \,  \ct_\emptyset) \\
 @A\phi_2AA    @AA\phi_2 A    @AA\phi_2 A \\
  (M_2(k)_{01}, \,  M_2(k)_b)  @>>>   (M_2(k)_{01}, \,  
M_2(k)_1)  @<<<   (M_2(k)_0 , \,  M_2(k)_\emptyset)    \, .
  \end{CD}
  \end{equation} 
The vertical maps are weak homotopy equivalences of pairs, by Theorem \ref{config}.  The horizontal  maps are induced by the values of the $\cf$-spaces on the morphisms in $\cf$. 
 
 For ease of notation, for a pair $(A,B)$ we write $A/B$ for the homotopy cofiber (mapping cone) $A \cup cB$.     By plugging in the values of these $\cf$-spaces, and taking homotopy cofibers, we get a commutative diagram
 \begin{equation}\label{downcofiber}
 \begin{CD}
 M_1 \times M_1 \times D^k/M_1 \times M_1 \times S^{k-1}  @>>>  M_1 \times M_1 \times D^k/ \fmak    @<\simeq << D(\tau_{M_1} \oplus \eps^k)/S(\tau_{M_1} \oplus \eps^k)  \\
 @V\phi_1 VV  @VV\phi_1 V  @VV\phi_1 V \\
 \ct_{01}/ \ct_b  @>>>    \ct_{01}/ \ct_1   @<<\simeq<    
\ct_0 / \ct_\emptyset  \\
 @A\phi_2AA    @AA\phi_2 A    @AA\phi_2 A \\
M_2 \times M_2 \times D^k/M_2 \times M_2 \times S^{k-1}  @>>>  M_2\times M_2 \times D^k/ \fmbk    @<\simeq << D(\tau_{M_2} \oplus \eps^k)/S(\tau_{M_2} \oplus \eps^k)   
\end{CD}
\end{equation}
The right hand horizontal maps are equivalences, because the commutative squares
defined by the $\cf$-spaces $\mokb$ and $\mtkb$ are pushouts, and therefore the commutative square defined by the $\cf$-space $\ct_\bullet$ is a homotopy pushout.
 By inverting these homotopy equivalences, as well as those induced by $\phi_2$, we get a homotopy commutative square,
 \begin{equation}
\begin{CD}\label{intdiag}
\Sigma^k((M_1 \times M_1)_+)  @>m_1 >> \Sigma^k(M_1 ^{\tau_{M_1}})\\
@V f_k  VV  @VV f_k V \\
\Sigma^k((M_2 \times M_2)_+)  @>m_2 >> \Sigma^k(M_2^{\tau_{M_2}})
\end{CD}
\end{equation}  Here the maps $f_k$ have the homotopy type of $\phi_2^{-1} \circ \phi_1$.  The right hand spaces are the suspensions of the  Thom spaces of the tangent bundles of $M_1$ and $M_2$ respectively.

Notice that property (2) in Theorem \ref{config} regarding the morphisms $\phi_1$ and $\phi_2$ and the mapping cylinder $\ct_{0,1}$ impies that the left hand map
$f_k : \Sigma^k((M_1 \times M_1)_+)  \to \Sigma^k((M_2 \times M_2)_+)  $ is given by the $k$-fold suspension of the equivalence $f\times f : M_1 \times M_1 \xr{\simeq} M_2 \times M_2$.  
Consider the right hand vertical equivalence,    $f_k: \Sigma^k(M_1)^{\tau_{M_1}} \to \Sigma^k(M_2^{\tau_{M_2}}).$  By diagram (\ref{downstairs}) $f_k$ is  induced by a map of pairs, 
$(D(\tau_{M_1} \oplus \eps^k), \, S(\tau_{M_1} \oplus \eps^k))  \to (D(\tau_{M_2} \oplus \eps^k), \, S(\tau_{M_2} \oplus \eps^k)) $  which on the ambient space,
$D(\tau_{M_1} \oplus \eps^k) \to D(\tau_{M_2} \oplus \eps^k)$ is homotopic to the map determined by
$f : M_1 \to M_2$ as in  property 3 of  Theorem \ref{config}.  Therefore the induced map in cohomology $(f_k)^* : h^*(\Sigma^k(M_2)^{\tau_{M_2}}) \cong  h^*(\Sigma^k(M_1)^{\tau_{M_1}})$
is an isomorphism  as modules over $h^*(M_2)$,  where the module structure on $h^*(\Sigma^k(M_1)^{\tau_{M_1}})$ is via  the isomorphism $f^* :h^*(M_2) \cong h^*(M_1)$.

 Moreover the isomorphism $(f_k)^* : h^*(\Sigma^k(M_2 ^{\tau_{M_2}})) \cong  h^*(\Sigma^k(M_1 ^{\tau_{M_1}}))$     preserves the Thom class in   cohomology.  To see this, notice that the
  horizontal maps  in  diagram (\ref{intdiag}) yield the intersection product in homology, after applying the Thom isomorphism. This implies
that the image of the fundamental classes $\Sigma^k([M_i] \times 1) \in h_{k+n}(\Sigma^k(M_i \times M_i))$ maps to the Thom classes
in $h_{k+n}(\Sigma^k(M_i)^{\tau_{M_i}}).$  Since the left hand vertical map is homotopic to $\Sigma^k (f \times f)$, and since the homotopy equivalence $f$ preserves the $h_*$-orientations, it preserves the fundamental classes.
Therefore by the commutativity of this diagram, $(f_k)_*$ preserves the Thom class.  These facts imply that after applying the Thom isomorphism, the isomorphism  $(f_k)^*$  is given by $f^* : h^*(M_2) \xr{\cong} h^*(M_1)$.

 This observation will be useful, as we will eventually   lift the map of  $\cf$-spaces given in Theorem \ref{config}    up to the level of loop spaces, and we'll consider the analogue of the diagram (\ref{intdiag}).   

To understand why this is relevant, recall 
  from \cite{chassullivan}, \cite{cohenjones} that the loop homology product 
$\mu : h_*(LM) \times h_*(LM) \to h_*(LM)$ can be defined in the following way.
Consider the pullback square
$$
\xymatrix{
\lmm \ar[r]^\iota \ar[d]_e & LM \times LM \ar[d]^{e \times e} \\
M \ar[r]_\Delta & M \times M
}
$$
where $e : LM \to M$ is the fibration given by evaluation at the basepoint: $e(\gamma) = \gamma (0)$.  Let $\eta (\Delta) $ be a tubular neighborhood
of the diagonal embedding of $M$, and let $\eta (\iota)$ be the inverse image
of this neighborhood in $LM \times LM$.  The normal bundle of $\Delta$ is the tangent bundle, $\tau_M$.  Recall that the evaluation map $e : LM \to M$ is a locally trivial fiber bundle \cite{klingenberg}.   Therefore   the tubular neighborhood $\eta (\iota)$ of $\iota : \lmm \hk LM \times LM$ is homeomorphic to total space of the pullback of the tangent bundle, $e^*(TM)$.  We therefore have a Pontrjagin-Thom collapse map,
\begin{equation}\label{tau}
\tau : LM \times LM \la LM \times LM /( (LM \times LM) - \eta (\iota)) \cong (\lmm)^{\tau_M}
\end{equation}
where $(\lmm)^{\tau_M}$ is the Thom space of the pullback $e^*(\tau_M) \to \lmm$.

 Now as pointed out in \cite{chassullivan}, there is a natural map 
$$
j\: \lmm \to LM
$$
given by sending a pair of loops $(\alpha, \beta)$ with the same starting point, to the concantenation
of the loops, $\alpha *\beta$.    
 The loop homology product is then defined to be the composition
\begin{equation}\label{mu}
\begin{CD}
 \mu_* \: h_*(LM \times LM) @>  \tau_* >>
 h_{*}((\lmm)^{TM}) @>\cap u> \cong> h_{*-n}(\lmm) @> j_* >>  h_{*-n}(LM) 
\end{CD} 
\end{equation} where $\cap u$ is   the Thom isomorphism given by capping with the Thom class.

\med
Now consider the fiber bundles, $LM_i \times LM_i \times D^k \to M_i \times M_i \times D^k = M_i(k)_{01}$ for $i = 1,2$.  
By restricting this bundle to the spaces 
$M_i(k)_j$, $j \in \text{\rm Ob}(\cf)$, we obtain   
$\cf$-spaces which we call $\cl M_i(k)_\bullet$, for $i = 1,2$. So we have morphisms of $\cf$-space, $e: \cl M_i(k)_\bullet \to  M_i(k)_\bullet$ which on every object is a fiber bundle, and every morphism induces a pull-back square.  

Similarly, let $\cl \ct_\bullet$ be the $\cf$ space obtained by restricting the fibration
$L(T_{f\times f}) \times D^k \xr{e \times 1}T_{f\times f} \times D^k = \ct_{01}$ 
to the spaces $\ct_j$, for $j \in \text{\rm Ob}(\cf)$.  

The morphisms  $\phi_i$ of Theorem \ref{config} then lift to give weak equivalences of $\cf$-spaces, $\cl \phi_i : \cl M_i(k)_\bullet  \to \cl \ct_\bullet$ that make the following diagram of $\cf$-spaces commute:
\begin{equation}\label{ell}
\begin{CD}
\cl \mokb  @>\cl\phi_1  >> \cl\ct_\bullet @<\cl\phi_2 << \cl \mtkb \\
@VeVV   @VVe V @VVeV \\
  \mokb  @> \phi_1  >>  \ct_\bullet @< \phi_2 <<   \mtkb 
  \end{CD}
  \end{equation}

The commutative diagram of maps of pairs (\ref{downstairs}) lifts to give a corresponding diagram
with spaces $\cl M_i(k)_\bullet$ replacing $M_i(k)_\bullet$, and $\cl \ct_\bullet$ replacing $\ct_\bullet$.  There is also a corresponding commutative diagram of quotients, that lifts
the diagram (\ref{downcofiber}).  The result is a homotopy commutative square, which lifts
square (\ref{intdiag}). 
\begin{equation}
\begin{CD}\label{liftint}
\Sigma^k((LM_1\times LM_1)_+)  @>\tau_1 >> \Sigma^k(LM_1 \times_{M_1} LM_1)^{\tau_{M_1}}) \\
@V \tilde f_k  VV  @VV \tilde f_k V \\
\Sigma^k((LM_2\times LM_2)_+)  @>\tau_1 >> \Sigma^k(LM_2 \times_{M_2} LM_2)^{\tau_{M_2}})
\end{CD}
\end{equation} Here the maps $\tilde f_k$ have the homotopy type of $L\phi_2^{-1} \circ L\phi_1$.   

Now as argued above, the description of the maps $L\phi_i : \cl M_i(k)_{(0,1)} \to \cl \ct_{(0,1)}$, that is, 
$LM_i \times LM_i \times D^k \to LT_{f\times f} \times D^k$ as the loop functor applied to the inclusion as the ends of the cylinder, implies that the map
$\tilde f_k : \Sigma^k((LM_1\times LM_1)_+)  \to  \Sigma^k((LM_2\times LM_2)_+) $ is homotopic to the $k$-fold suspension of $Lf \times Lf : LM_1 \times LM_1 \xr{\simeq} LM_2 \times LM_2.$ Moreover, in cohomology, the map $\tilde f_k^* : h^*(\Sigma^k(LM_2 \times_{M_2} LM_2)^{\tau_{M_2}}) \to h^*(\Sigma^k(LM_1 \times_{M_1} LM_1)^{\tau_{M_1}})$ preserves   Thom classes because the bundles are pulled back from
bundles over $M_1$ and $M_2$ respectively, and as seen above, $(f_k)^* : h^*(\Sigma^k(M_2 ^{\tau_{M_2}})) \cong  h^*(\Sigma^k(M_1 ^{\tau_{M_1}}))$     preserves
Thom classes.  Also, since this map is, up to homotopy,  induced by a map of pairs
$$
 L\phi_2^{-1} \circ L\phi_1: (D(e^*(\tau_{M_1} \oplus \eps^k)), \, S(e^*(\tau_{M_1} \oplus \eps^k)))  \to (D(e^*(\tau_{M_2} \oplus \eps^k)), \, S(e^*(\tau_{M_2} \oplus \eps^k))) 
 $$
   it  induces an isomorphism of $h^*(D(e^*(\tau_{M_2} \oplus \eps^k)))$ modules, where this ring acts on $ h^*(\Sigma^k(LM_1 \times_{M_1} LM_1)^{\tau_{M_1}})$ via the homomorphism, $ (L\phi_2^{-1} \circ L\phi_1)^* : h^*(D(e^*(\tau_{M_2} \oplus \eps^k))) \xr{\cong} h^*(D(e^*(\tau_{M_1} \oplus \eps^k)))$.  But by the lifting of property 3 in Theorem \ref{config}, this map is homotopic to the compositon,
 $$
 (D(e^*(\tau_{M_1} \oplus \eps^k)) \xr{\rm project} LM_1 \times_{M_1}LM_1 \xr{Lf \times Lf} LM_2 \times_{M_2}LM_2 \xr{\rm zero}(D(e^*(\tau_{M_2} \oplus \eps^k)).
 $$
 Hence when one applies the Thom isomorphism to both sides, the isomorphism 
 $$\tilde f_k^* : h^*(\Sigma^k(LM_2 \times_{M_2} LM_2)^{\tau_{M_2}}) \to h^*(\Sigma^k(LM_1 \times_{M_1} LM_1)^{\tau_{M_1}})$$
 is given by $(Lf \times Lf)^* : h^*(LM_2 \times_{M_2}LM_2) \xr{\cong} h^*(LM_1 \times_{M_1}LM_1)$.

 \med
By the definition of the loop product (\ref{mu}), to prove that $(Lf)_* : h_*(LM_1) \to h_*(LM_2)$ is a ring isomorphism,  we need to show that the diagram
$$
\begin{CD}
h_*(LM_1\times LM_1)  @>(\tau_1)_* >>  h_{*}((LM_1 \times LM_1)^{\tau_{M_1}}) @>\cap u> \cong> h_{*-n}((LM_1 \times_{M_1} LM_1)) @> j_* >>  h_{*-n}(LM_1)  \\
@V(Lf \times Lf)_* VV @V(\tilde f_k)_* VV     @VV(Lf \times Lf)_* V  @VV(Lf)_* V \\
h_*(LM_2\times LM_2)  @>(\tau_2)_*>>  h_{*}((LM_2 \times LM_2)^{\tau_{M_2}}) @>\cap u> \cong> h_{*-n}((LM_2 \times_{M_2} LM_2)) @> j_* >>  h_{*-n}(LM_2)
 \end{CD}
$$
 commutes.  We have now verified that the left and middle squares commute.  But the right hand square obviously commutes.  Thus $(Lf)_* : h_*(LM_1) \to h_*(LM_2)$ is a ring isomorphism as claimed.   
 
 To prove that $Lf$ is a map of $BV$- algebras, recall that the $BV$-operator $\Delta$ is defined in terms of the $S^1$-action.   Clearly $Lf$ preserves this action, and hence induces an isomorphism of $BV$-algebras.  This will imply that 
  $Lf$ induces an isomorphism of the string Lie algebras for the following reason.  Recall the definition of the Lie bracket from \cite{chassullivan}.   Given $\alpha \in h_p^{S^1}(LM)$ and $\beta \in h_q^{S^1}(LM)$, then the bracket $[\alpha, \beta]$ is the image of $\alpha \times \beta$ under the composition,
 
\begin{align}
h_p^{S^1}(LM) \times h_q^{S^1}(LM) &\xr{\text{\rm tr}_{S^1} \times  \text{\rm tr}_{S^1}} h_{p+1}(LM) \times h_{q+1}(LM)\notag \\
&\xr{\rm loop\,  product} 
h_{p+q+2-n}(LM) \xr{j} h_{p+q+2-n}^{S^1}(LM). 
\end{align}
Here $\text{\rm tr}_{S^1} \: h_*^{S^1}(LM) \to h_{*+1}(LM)$ is the $S^1$ transfer map 
(called ``M" in \cite{chassullivan}), 
and $j \: h_*(LM) \to h^{S^1}_*(LM)$ is the usual 
map that descends nonequivariant homology to 
equivariant homology (called ``E" in \cite{chassullivan}). 
We refer the reader to \cite{Adem-Cohen-Dwyer} for a concise definition of the $S^1$-transfer.

We now know that $Lf$ preserves the loop product, and since it is an $S^1$-equivariant map, it preserves the transfer map $\text{\rm tr}_{S^1}$ and the map $j$.  Therefore it preserves the string bracket.  
 \end{proof}

\med
\section{Relative embeddings and the proof of Theorem \ref{main}}

Theorem \ref{config} reduces the proof of the
homotopy invariance of the loop   product and the string bracket 
(Theorem \ref{main})
to the homotopy invariance of the 
${\cal F}$-spaces associated with the
embeddings diagonal of $M_1$ and $M_2$.
The goal of the present section is to prove 
Theorem \ref{config}.

\subsection{Relative smooth embeddings}

Let $N$ be a compact smooth manifold of dimension $n$
whose boundary $\partial N$ comes equipped with a smooth manifold decomposition
$$
\partial N = \partial_0 N \cup \partial_1 N
$$
in which $\partial_0 N$ and $\partial_1 N$ are glued together
along their common boundary 
$$
\partial_{01} N \,\, := \,\, \partial_0 N \cap \partial_1 N\, .
$$
Assume that $K$ is a space obtained
from $\partial_0 N$ by attaching a finite number of cells.
Hence we have a relative cellular complex
$$
(K,\partial_0 N)\, . 
$$
It then makes sense to speak of the {\it relative dimension}
$$
\dim (K,\partial_0 N) \le k
$$ as being the maximum dimension of the
attached cells.

Let $$f\: K \to N$$ be a map of spaces which extends the identity map of $\partial_0 N$.

\begin{definition} We call these data, $(K, \p_0N, f \: K\to N)$, a {\it relative smooth embedding problem}
\end{definition}

\begin{definition} 
A {\it solution} to the relative smooth embedding problem consists
of  
\begin{itemize}
\item a codimension zero compact submanifold
$$
W \subset N
$$
such that $\partial W \cap \partial N = \partial_0 N$ and this intersection
is transversal, and 
\item a homotopy of $f$, fixed
on $\partial_0 N$, to a map of the form
$$
\begin{CD}
K @>\sim >> W @>\subset >> N
\end{CD}
$$
in which the first map is a homotopy equivalence.
\end{itemize}
\end{definition}

\begin{lemma} \label{smooth} 
If $2k < n$, then there is a solution to the relative smooth
embedding problem.
\end{lemma}


We remark that Lemma 4 is essentially a simplified version of a result
of Hodgson \cite{Hodgson} who strengthens it by $r$ dimensions when 
the map $f$ is $r$-connected.

\begin{proof}[Proof of the Lemma \ref{smooth}]  
First assume that $K = \partial_0 N \cup D^k$ is the effect
of attaching a single $k$-cell to $\partial_0 N$. Then the restriction of
$f$ to the disk gives a map
$$
(D^k,S^{k-1}) \to (N,\partial_0 N)
$$
and, by transversality, we
can assume that its restriction $S^{k-1} \to \partial_0 N$ 
is a smooth embedding. 
Applying transversality again, the map on $D^k$ 
can be generically deformed relative to
$S^{k-1}$ to a smooth embedding. Call the resulting embedding $g$.
Let $W$ be defined by taking a regular neighborhood of 
$\partial_0 N \cup g(D^k) \subset N$. Then $g$ and $W$ give 
the desired solution
in this particular case.

The general case is by induction on the 
set of cells attached to $\partial_0 N$. 
The point is that if a solution $W \subset N$ 
has already been achieved on a subcomplex $L$ of
$K$ given by deleting one of the top cells, then  removing  the interior of
$W$ from $N$  gives a new manifold $N'$, such that $\partial N'$ 
has a boundary decomposition. The
attaching map $S^{k-1} \to L$ can be deformed (again using transversality)
to a map into $\partial_0 N'$. Then we have reduced
to a situation of solving the problem for a map of
the form $D^k \cup \partial_0 N' \to N'$, which we know can be solved by the previous paragraph.
\end{proof}

\med
We now thicken the complex $K$ by crossing $\p_0N$ with a disk.  Namely, 
for an integer $j \ge 0$,  define the space
$$
K_j \simeq K \cup_{\partial_0 N} (\partial_0 N) \times D^j\, ,
$$
where we use the inclusion $\partial_0 N \times 0 \subset
(\partial_0 N) \times D^j$ to form the amalgamated union.
Then $(K,\partial_0 N) \subset (K_j,(\partial_0 N) \times D^j)$ is
a deformation retract, and the map $f\:K \to N$ extends in the evident
way to  a map
$$
f_j \: K_j \to N \times D^j
$$
that is fixed on $(\partial_0 N) \times D^j$.

\begin{theorem} \label{smooth+j} 
Let  $f\: K \to N$ be as above, but without the dimension
restrictions.  Then for sufficiently large    $j \ge 0$, the   
  embedding problem for the map $f_j\: K_j \to N \times D^j$ admits a solution.
\end{theorem}

\begin{proof} The relative dimension of $(K_j, \, (\partial_0 N) \times D^j)$
is $k$, but for sufficiently large $j$ we  have $2k \le n + j$.  The result follows from the previous lemma.
\end{proof}

\subsection{Relative Poincar\'e embeddings}
Now suppose more generally that $(N,\partial N)$ 
is a (finite) Poincar\'e pair of dimension $n$ equipped with
a {\it boundary decomposition}
such that $\partial_0 N$ is a smooth manifold.  By this, we mean
we have an expression of the form
$$
\partial N \,\, = \,\, \partial_0 N \cup_{\partial_{01} N} \partial_1 N
$$
in which $\partial_0 N$ is a manifold with boundary $\partial_{01}N$ and
also $(\partial_1 N,\partial_{01}N)$ is a Poincar\'e pair. Furthermore,
we assume that the fundamental classes for each of theses pairs glue
to a fundamental class for $\partial N$.  These fundamental classes lie in ordinary homology.

As above, let 
$$
f\: K \to N
$$
be a map which is fixed on $\partial_0 N$. We will assume
that the relative dimension of $(K,\partial_0 N)$ is 
at most $n-3$. Call these data a 
{\it relative Poincar\'e embedding problem}.

\begin{definition} A {\it solution}  of a relative Poincar\'e embedding problem as above
consists of
\begin{itemize} 
\item a Poincare pair $(W,\partial W)$, and a Poincar\'e decomposition
$$
\partial W = \partial_0 N \cup_{\partial_{01} N} \partial_1 W
$$
such that  each of the maps $\partial_{01} N \to \partial_1 W$ and
$\partial_0 N  \to W$ is $2$-connected;
\item a Poincare pair $(C,\partial C)$ with Poincar\'e
decomposition
$$
\partial C = \partial_1 W \cup_{\partial_{01} N} \partial_1 N \, ;
$$
\item a weak equivalence $h\: K \to W$ which is fixed on $\partial_0 N$;
\item a weak equivalence
$$
e\:  W \cup_{\partial_1 W} C \to N
$$
which is fixed on $\partial N$, such that
$e\circ h$ is homotopic to $f$ by a homotopy fixing $\partial_0 N$.
\end{itemize}
\end{definition}

The above is depicted in the following schematic
homotopy decomposition of $N$:
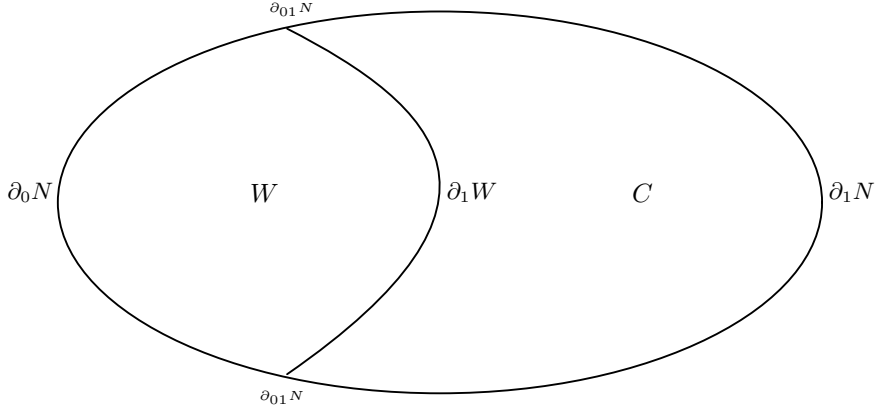
\bigskip

$$
\begin{texdraw}
\lellip rx:2 ry:1
\htext (-2.27 0) {\small $\partial_0 N$}
\htext (-1 0) {$W$}
\htext (.03 0) {\small $\partial_1 W$}
\htext (1 0) {$C$}
\htext (2.03 0) {\small $\partial_1 N$}
\htext (-.9 .97) {\tiny $\partial_{01}N$}
\htext (-.95 -1.06) {\tiny $\partial_{01}N$}
         \move (-.8 .91)
         \clvec (0 .5)(.5 0)(-.8 -.9)
\end{texdraw}
$$
\bigskip

The space $C$ is called the {\it complement},
which is a Poincar\'e space with boundary $\partial_1 W \cup \partial_1 N$.
The above spaces 
assemble to give a  strictly commutative square which is  
homotopy cocartesian:  
\begin{equation} \label{the_diagram}
\xymatrix{
(\partial_1 W,\partial_{01} N ) \ar[r] \ar[d] & (C,\partial_1 N) \ar[d] \\
(W,\partial_0 N) \ar[r]   & (N,\partial N)
}
\end{equation}
(compare \cite{Klein_haef2}).
From here through the rest of the paper we refer to 
such a commutative square as a ``homotopy pushout".
\med

As above, we can construct maps $f_j \: K_j \to N \times D^j$, which  define a family of relative Poincare embedding problems.  Our goal in this section is to prove the analogue of Theorem \ref{smooth+j} that shows that for sufficiently large $j$ one can find solutions to these problems.  

\med
We begin with the following result, comparing the smooth to the Poincare relative embedding problems.

\begin{lemma} \label{smooth=>poincare} Assume that is $M$ is a compact
smooth manifold equipped with a boundary decomposition. Let
$$
\phi\: (N;\partial_0 N,\partial_1 N) \to (M,\partial_0 M,\partial_1 M)
$$
be a homotopy equivalence whose restriction
$\partial_0 N \to \partial_0 M$ is a diffeomorphism.

Then the relative Poincar\'e embedding problem for $f$ admits
a solution if
the relative smooth embedding problem for $\phi\circ f$ admits a solution.
\end{lemma}

\begin{proof} A solution of the smooth problem 
together with a choice of homotopy inverse for $h$ extending
the inverse diffeomorphism on $\partial_0 M$ gives  solution
to the relative Poincar\'e embedding problem.
\end{proof}

\med
Now suppose that
$D(\nu) \to \partial_0 N$ is the unit disk bundle of the 
normal bundle of an embedding of $\p N$ into codimension $\ell$ Euclidean space, $\br^{n + \ell}$. 
The zero section then gives  an inclusion
$\partial_0 N \subset D(\nu)$.
Set 
$$
K_\nu := K \cup_{\partial_0 N} D(\nu).
$$
Clearly, $K_\nu$ is canonically homotopy equivalent to $K$.

Assuming $\ell$ is sufficiently large, there exists   
Spivak normal fibration \cite{Spivak}
$$
S(\xi) \to N
$$  
whose fibers have the homotopy type of an $\ell-1$ dimensional sphere.  
Then,  by the uniqueness of the 
Spivak fibration \cite{Spivak}, 
we have  a fiber homotopy equivalence over $\partial N$
$$
\begin{CD}
S(\nu) @> \simeq >> S(\xi_{|\partial_0 N}).
\end{CD}
$$
Let $D(\xi)$ denote the fiberwise cone fibration of $S(\xi) \to N$. Then we
have a canonical  map
$$
f_\nu\: K_\nu \to D(\xi)
$$
which is fixed on the $D(\nu)$. Note that
$$
\partial D(\xi) = D(\nu) \cup S(\xi)
$$
is a decomposition of Poincar\'e spaces such that
$D(\nu)$ has the structure of a smooth manifold.
Let us set $\partial_0 D(\xi) := D(\nu)$  and
$\partial_1 D(\xi) := S(\xi)$.

Then the classical construction of the Spivak fibration (using 
regular neigbhorhood theory in Euclidean space)
shows that there is a homotopy equivalence 
 $$
\begin{CD}
(D(\xi);\partial_0 D(\xi),\partial_1 D(\xi)) 
@> \sim >> (M;\partial_0 M,\partial_1 M)
\end{CD}
$$
in which $M$ is a compact codimension zero submanifold of
some Euclidean space. Furthermore, the restriction
$\partial_0 D(\xi) \to \partial_0 M$ is a diffeomorphism.

Consequently, by Lemma \ref{smooth} and lemma \ref{smooth=>poincare} we obtain

\begin{proposition} \label{f-nu} If the rank of $\nu$ is sufficiently large, then
the relative Poincar\'e embedding problem for
$f_\nu$ has a solution.
\end{proposition}

Let $\eta$ denote a choice  of inverse
for $\xi$ in the Grothendieck group of reduced spherical fibrations over $N$.  
For simplicity, we may assume that the fiber of 
$\eta$ is a sphere of dimension $\dim N - 1$. Then $\xi$ restricted to
$\partial_0 N$ is fiber homotopy equivalent to
$\tau_{\partial_0 N} \oplus \epsilon$, where $\tau_{\partial_0 N}$ is
the tangent sphere bundle of $\partial_0 N$ and $\epsilon$ is the
trivial bundle with fiber $S^0$. For simplicity, we will
assume that $\xi$ restricted to $\partial_0 N$ has been identified
with $\tau_{\partial_0 N} \oplus \epsilon$. Similarly, we
will choose an identification of $\xi_{\partial_1 N}$ with
$\tau_{\partial_1 N} \oplus \epsilon$, where $\tau_{\partial_1 N}$
is any spherical fibration over $\partial_1 N$ that represents 
the Spivak tangent fibration.

Since $\xi \oplus \eta $ is trivializable, 
for some integer $j$ we get a homotopy equivalence
$$
(D(\xi \oplus \eta); D(\nu \oplus \tau_{\partial_0 N});
D(\nu\oplus \tau_{\partial_1 N}) \cup S(\xi\oplus \eta)) \to
(N\times D^j; (\partial_0 N) \times D^j;(\partial_1 N) \times D^j \cup 
N \times S^{j-1}) 
$$
which restricts to a diffeomorphism
  $D(\nu\oplus \tau_{\partial_0 N}) \to (\partial_0 N) \times D^j$.

Now a choice of solution of the relative Poincar\'e embedding problem
for $f_\nu$, as given by Proposition \ref{f-nu}, guarantees that 
the relative problem for $f_{\nu\oplus \tau}$ has a solution.
But clearly, the latter is identified with the map
$f_j\: K_j \to N \times D^j$. Consequently, we have proven the following.

\begin{theorem} \label{solution} If $j \gg 0$ is sufficiently large, then the
relative Poincar\'e embedding problem for $ f_j\: K_j \to N \times D^j$ has a solution.
\end{theorem}

\subsection{Application to diagonal maps and a proof of Theorem \ref{config}}

We now give a proof of Theorem \ref{config}.  By the results of section 1, this will complete the proof of Theorem \ref{main}.

\med
Let $f \: M_1 \to M_2$ be a homotopy equivalence of
  closed smooth manifolds.
Using an identification of the tangent bundle $\tau_{M_1} $ with the normal bundle of the diagonal, $\Delta \: M_1 \to M_1 \times M_1$,   we have an embedding
$$
D(\tau_{M_1}) \subset M_1\times M_1\, ,
$$
which is identified with a compact tubular neighborhood of
the diagonal.   
The closure of its complement will be denoted $F(M_1, 2)$.   Notice  that
the inclusion $F(M_1, 2) \subset M_1^{\times 2} - \Delta$ is a weak
equivalence of spaces over $M_1^{\times 2}$ (i.e., it is a morphism
of spaces over $M_1^{\times 2}$ whose underlying map of spaces
is a weak homotopy equivalence). Notice also  that we have a decomposition
$$
M_1^{\times 2} =  D(\tau_{M_1}) \cup_{S(\tau_{M_1})} F(M_1, 2).
$$
Making the same construction  with $M_2$, we also have a decomposition
$$
M_2^{\times 2} =  D(\tau_{M_2}) \cup_{S(\tau_{M_2})} F(M_2, 2) \, .
$$

Notice that since $f \: M_1 \to M_2$ is a homotopy equivalence, 
  the composite 
$$
\begin{CD}
D(\tau_{M_1}) @> \text{projection} >>  M_1 @> f >> M_2 @> 
\text{zero section} >\hk> D(\tau_{M_2})
\end{CD}
$$
is also a homotopy equivalence.
Let $T$ be the {\it mapping cylinder} of this composite  map. Then we
have a pair
$$
(T,D(\tau_{M_1}) \amalg D(\tau_{M_2}))\, .
$$
Furthermore, up to homotopy, we have a preferred identification of $T$ with the mapping cylinder of $f$.

The map $f^{\times 2} \: M_1^{\times 2} \to M_2^{\times 2}$ also
has a mapping cylinder $T^{(2)}$ which contains the manifold
$$
\partial T^{(2)} := M_1^{\times 2} \amalg M_2^{\times 2}.
$$
Then  $(T^{(2)},\partial T^{(2)})$ is a Poincar\'e pair.
Furthermore, 
$$
\partial T^{(2)} = (D(\tau_{M_1}) \amalg D(\tau_{M_2})) \cup (F(M_1, 2) \amalg F(M_2, 2))
$$
is a manifold decomposition. Let us set $\partial_0 T^{(2)} = 
D(\tau_{M_1}) \amalg D(\tau_{M_2})$ and $\partial_1 T^{(2)} = (F(M_1, 2) \amalg F(M_2, 2))$.

Since the diagram
$$
\xymatrix{
 M_1 \ar[r]^f \ar[d]_{\Delta} & M_2
\ar[d]^{\Delta} \\
 M_1 \times M_1 \ar[r]_{f\times f} & M_2 \times M_2 
 }
$$
commutes, we get an induced map of mapping cylinders. 
This map, together with our preferred identification of the 
cylinder of $f$ with
$T$,  allows the construction of a map
$$
g\:T \to T^{(2)}
$$
which extends the identity map of $\partial_0 T^{(2)}$. 
In other words, $g$ is a relative Poincar\'e embedding problem.

By Proposition \ref{f-nu}, there exists an integer $j \gg 0$ such that
the associated relative Poincar\'e embedding problem
$$
g_j\: T_j \to T^{(2)} \times D^j
$$
has a solution. Here,
$$
T_j \,\, := \,\, T \cup \, (\partial_0 T^{(2)}) \times D^j\, ,
$$  
and 
$$
\p_0(T^{(2)} \times D^j) := \left(D(\tau_{M_1}) \times D^j \right) \amalg  \left(D(\tau_{M_2}) \times D^j \right) \qquad 
\p_1(T^{(2)} \times D^j) :=  F_{D^j}(M_1, \, 2) \amalg F_{D^j}(M_2, \,2).
$$
where, for convenience, we are redefining $F_{D^j}(M, \, 2)$
as $M \times M \times S^{j-1} \cup F(M,2) \times D^j$ (cf.\ S1).

This makes $T^{(2)}\times D^j$ a Poincar\'e space with boundary decomposition
$$
\partial (T^{(2)} \times D^j)\,\,  = \partial_0 (T^{(2)} \times D^j) \cup
\partial_1 (T^{(2)} \times D^j)\, .
$$

\med
By definition \ref{solution}, a solution to this Poincar\'e embedding problem  yields  Poincar\'e  pairs  $(W, \p W)$ and $(C, \p C)$,  with the following properties. 

\begin{itemize}\label{solved}
\item
 $
\p W = \p_0(T^{(2)}\times D^j) \cup \p_1W$, where $\p_0 ({T^{(2)}\times D^j}) = \left(D(\tau_{M_1}) \times D^j \right) \amalg  \left(D(\tau_{M_2}) \times D^j \right)$ and $\p_1 W \hk W$ is 2-connected,
\item  $\p C = \p_1 W \cup \p_1(T^{(2)}\times D^j), $ where $\p_1(T^{(2)}\times D^j) =  F_{D^j}(M_1, \, 2) \amalg F_{D^j}(M_2, \,2).$  Notice that $\p_{01}(T^{(2)}\times D^j) = \p(D(\tau_{M_1} \times D^j)) \, \amalg \, \p(D(\tau_{M_2} \times D^j)).$
\item
There is a  weak equivalence, $h \: T_j \xr{\simeq} W$, fixed on $\left(D(\tau_{M_1}) \times D^j \right) \amalg  \left(D(\tau_{M_2}) \times D^j \right)$.   
 \item There is  a weak equivalence
$$
e \: W \cup_{\p_1W}C \to T^{(2)}\times D^j
$$
which is fixed on $\p(T^{(2)}\times D^j)$,   such that $e \circ h$ is homotopic to 
$g_j \: T_j \to T^{(2)}\times D^j$ by a homotopy fixing $\left(D(\tau_{M_1}) \times D^j \right) \amalg  \left(D(\tau_{M_2}) \times D^j \right)$.
\end{itemize}
The above homotopy decomposition of $T^{(2)}\times D^j$ is indicated in the
following schematic diagram:
\medskip

$$
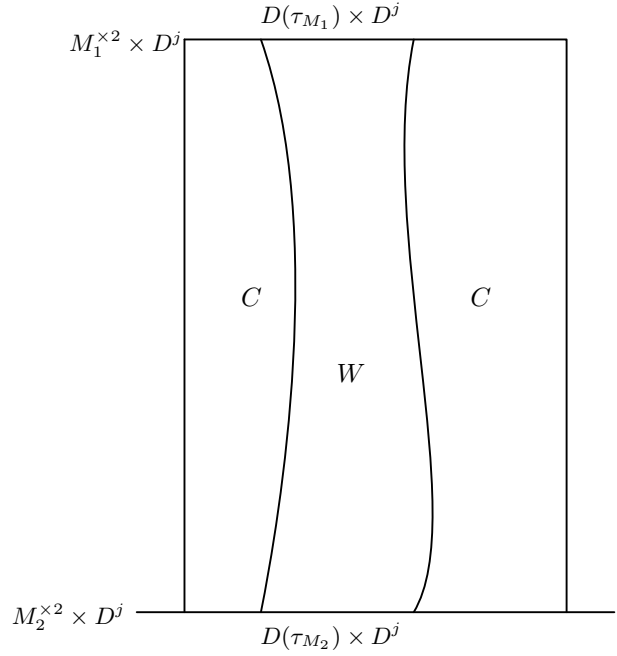
\btexdraw
\lvec (2 0)
\lvec (2 -3)
\move (2.25 -3)
\lvec (-.25 -3) 
\move (0 -3)
\lvec (0 0)
\move (.4 0)
\clvec (.75 -1)(.5 -2.5)(.4 -3)
\move (1.2 0)
\clvec (1 -1)(1.5 -2.5)(1.2 -3) 
\htext (.8 -1.8) {$W$}
\htext (.3 -1.4) {$C$}
\htext (1.5 -1.4) {$C$}
\htext (.4 .05) {\small $D(\tau_{M_1}) \times D^j$}
\htext (.4 -3.2) {\small $D(\tau_{M_2}) \times D^j$}
\htext (-.6 -.1) {\small $M_1^{\times 2} \times D^j$}
\htext (-.9 -3.1) {\small $M_2^{\times 2} \times D^j$}
\etexdraw
$$
\medskip
Furthermore, the complement $C$ and the
normal data $\partial_1 W$
of the solution sits in a commutative diagram
of pairs
$$
\xymatrix{
(M_1^{\times 2} \times S^{j-1},\emptyset) \ar[d]_\cap \ar[rr]^\subset_\sim 
&& (T^{(2)}\times S^{j-1},\emptyset) \ar[d]^\cap && 
(M_2^{\times 2}\times  S^{j-1},\emptyset) \ar[d]^\cap  \ar[ll]_\supset^\sim \\
(F_{D^j}(M_1,2),S(\tau_{M_1} + \epsilon^j))
 \ar[rr]^\subset  \ar[d]_{\cap} && (C,\partial_1 W)  \ar[d]^e &&
(F_{D^j} (M_2, 2),S(\tau_{M_2} + \epsilon^j))  \ar[d]^{\cap}
\ar[ll]_{\supset}\\ 
(M_1^{\times 2} \times D^j,D(\tau_{M_1}) \times D^j) 
\ar[rr]^{\subset}_\sim && (T^{(2)} \times D^j,T_j)  
&&
(M_2^{\times 2}\times D^j,D(\tau_{M_2}) \times D^j) \ar[ll]_{\supset}^\sim \, .
}
$$
 Here  each (horizontal) 
arrow marked with $\sim$ is a weak homotopy equivalence.  
Each column describes an ${\cal F}$-space (cf.\ 
Definition \ref{fivespace}). In fact, the outer columns
are the ${\cal F}$-spaces $M_i(j)$  described in \S1.  
Furthermore, the horizontal maps describe morphisms of ${\cal F}$-spaces.

Consequently, to complete
the proof of Theorem \ref{config},  
it suffices to show these morphisms of ${\cal F}$-spaces
are weak equivalences.
We are therefore reduced to showing that the horizontal arrows in the
second row are weak homotopy equivalences.

By symmetry, it will suffice
to prove that the left map in the second row, 
$$(F_{D^j}(M_1,2),S(\tau_{M_1} + \epsilon^j)) \to (C,\partial_1 W)
$$   
is a weak  equivalence.
\med 

We will prove that the map $F_{D^j}(M_1,2)\to C$
is a weak equivalence; the proof that $S(\tau_{M_1} + \epsilon^j) \to 
\partial_1 W$ is a weak equivalence is similar and will be
 left to the reader. 
To do this, consider the following commutative diagram.
\begin{equation}\label{bigcommute}
\xymatrix{
F_{D^j}(M_1, 2)  \ar[r]^= \ar[d]_{\cap} &F_{D^j}(M_1, 2) 
\ar[d]^\cap \ar[r]^{\hk} & C \ar[d]^e\\
M_1^{\times 2} \times D^j \ar[r]_>>>>{\hk}   & 
( M_1^{\times 2} \times D^j) \cup_{ (D(\tau_{M_1}) \times D^j)} W   \ar[r]_>>>>{\hk}  
& T^{(2)} \times D^j\, .
}
\end{equation}

 \med
 \begin{lemma}\label{push}
 Each of the commutative squares in diagram 
\eqref{bigcommute} is a homotopy pushout.
 \end{lemma}
 
 \med
 Before we prove this lemma, we show how we will use it to complete the proof of Theorem \ref{config}.
 
\begin{proof}[Proof  of Theorem \ref{config}]  \rm   
By the lemma, since each of the squares
 of this diagram is a homotopy pushout, then so is the outer diagram,
 $$
 \xymatrix{
 F_{D^j}(M_1, 2) \ar[r]^{\hk} \ar[d]_\cap &  C \ar[d]^e \\
 M_1^{\times 2} \times D^j  \ar[r]_\hk & T^{(2)} \times D^j.
 }
 $$
 Now recall that $ T^{(2)}$  is the mapping cylinder of the homotopy equivalence,
 $f^{\times 2}\: M_1^{\times 2} \to M_2^{\times 2}.$   Therefore the inclusion,
 $M_1^{\times 2}  \to T^{(2)}$ is an equivalence, 
and hence so is the bottom horizontal
 map in this pushout diagram, $M_1^{\times 2} \times D^j  \hk  T^{(2)} \times D^j$.
 Furthermore, the inclusion $ F_{D^j}(M_1, 2) \to M_1^{\times 2} \times D^j $ is $2$-connected, assuming the  dimension of $M$ is $2$ or larger.    Therefore by the pushout property of this square and the Blakers-Massey theorem,  we conclude that the top horizontal map
 in this diagram, $  F_{D^j}(M_1, 2)  \hk  C$ is a homotopy equivalence.

 As described before, this is what was needed to complete the  proof of  
Theorem \ref{config}.  \end{proof}

 \begin{proof}[Proof of Lemma \ref{push}]
 We first consider the right hand commutative square.  By the properties of the solution
 to the relative embedding problem given above in (\ref{solved}), we know that $e\: C \to T^{(2)}\times D^j$ extends to an equivalence, $e\: C \cup_{\p_1W}W \xr{\simeq}  T^{(2)}\times D^j$.  Now notice that the intersection of $\p_1 W$ with $F_{D_j}(M_1,2)$ is the boundary, $\p(D(\tau_{M_1}) \times D^j)$.  But 
 $$
 (F_{D_j}(M_1,2)) \cup_{\p(D(\tau_{M_1}) \times D^j)} W = ( M_1^{\times 2} \times D^j) \cup_{ (D(\tau_{M_1}) \times D^j)} W .
 $$
 This proves that the right hand square is a homotopy pushout.  
 
 We now consider the left hand diagram.  Again, by using the properties of the solution of  the relative embedding problem given above in (\ref{solved}), we know that
 the homotopy equivalence $h \: T_j \xr{\simeq} W$ extends to a homotopy equivalence,
 $$
 h  \: ( M_1^{\times 2} \times D^j) \cup_{ (D(\tau_{M_1}) \times D^j)} T_j  \xr{\simeq} ( M_1^{\times 2} \times D^j) \cup_{ (D(\tau_{M_1}) \times D^j)} W.
 $$
 But by construction, $T_j$ is homotopy equivalent to the mapping cylinder of the composite homotopy equivalence, $D(\tau_{M_1}) 
\xr{\text{\rm project}} M_1 \xr{f} M_2 \xr{\text{zero section}} D(\tau_{M_2})$.   This implies that the inclusion $D(\tau_{M_1}) \times D^j \hk T_j$ is a homotopy equivalence, and so $ ( M_1^{\times 2} \times D^j) \cup_{ (D(\tau_{M_1}) \times D^j)} T_j $ is homotopy equivalent to $ M_1^{\times 2} \times D^j$.  Thus the inclusion given by  the bottom horizontal map in the square in question, $M_1^{\times 2} \times D^j \hk (M_1^{\times 2} \times D^j) \cup_{ (D(\tau_{M_1}) \times D^j)} W$ is also a homotopy equivalence.  Since the top horizontal map is the identity,
 this square is also a homotopy pushout.  This completes the proof of Lemma \ref{push}, which was the last step in the proof of Theorem \ref{config}.  \end{proof}

\end{document}